\documentclass[english]{article}
\usepackage[T1]{fontenc}
\usepackage[latin9]{inputenc}
\usepackage{units}
\usepackage{amsmath}
\usepackage{amssymb}

\makeatletter
\usepackage{amsfonts}

\makeatother

\usepackage{babel}

\begin{document}

\title{Is Zero a Natural Number?}

\author{Peter Harremoës\\
Copenhagen Business College}

\maketitle

\section{Introduction}

The simplest way to answer the question in the title is to say that
zero by definition is a natural number, or, as other do, to say that
zero by definition is not a natural number. In mathematics we are
allowed to make the definitions as we like so it is just a matter
of definition and there is noting more to say about it.

If I agreed on that this would be a very short article, but if we
think of mathematics as a a set of mental and conceptual tools to
solve real world problems, some definitions are more useful than others
and then it make sense to discuss the consequences of different definitions,
and which definitions are the most useful ones. The most basic notions
of numbers are the cardinal and the ordinal numbers. Whether zero
should be considered as a natural number depends on whether we would
like to identify the natural numbers as the finite cardinal numbers
or with the finite ordinal numbers. Let us take a closer look at these
two classes of numbers.

\section{Cardinal numbers}

Many languages has several classes of numbers with cardinal numbers
as one of them. Although counting has been an important activity for
thousands of years and is one of the most basic models in mathematics,
it was not before Cantor that a well-founded mathematical understanding
emerged. The basic definition is that two sets are equipotent if there
is a one-to-one pairing of the elements. It is easy to see that we
get an equivalence relation. It is often inconvenient to use expressions
like \textquotedbl{}I have found some mushrooms and the set of mushrooms
is equipotent to the set of poles in our tent\textquotedbl{} because
maybe the one that hears this message does not know the tent that
is referred to that well. Although it is often inconvenient, the method
is still in use. For instance: make one prayer for each pearl on the
rosary. This makes sense without any counting as long as you have
a rosary. The next step is to use some standard set to refer to. This
is often related to the fingers and perhaps also the toes, because
this set is normally available and these standard sets are equipotent
for most humans. Most people have some more or less developed standard
sets that they use to understand the higher numbers. For instance
300 is for me is approximately the number of pupils in a Danish gymnasium.
In books on science one also finds expressions like \textquotedbl{}This
is of the same order as the number of atoms in the universe\textquotedbl{}.

The next step is to have standard labels for sets of each equipotency
class. Thus \textquotedbl{}five\textquotedbl{} is the label we use
for all sets with five elements. From this point of view it does not
make much sense to ask what \textquotedbl{}five\textquotedbl{} is.
It is just a word that we use when we want to tell that a set has
as many elements as we have fingers on one hand.

Most modern languages are far from fully developed in the use of such
standard labels. For instance in English one can either say \textquotedbl{}one
man\textquotedbl{} or \textquotedbl{}a man\textquotedbl{} where both
indicate the same number, but \textquotedbl{}a\textquotedbl{} is normally
not considered as a numeral. Some languages distinguish between different
classes of objects and use different numerals for different classes.
In some languages the nouns have different genders and the numeral
change depending on the gender. In many languages the numerals interact
with the inflections system of the language. Many languages also have
inflections in singular or plural. I some languages there is even
an inflection called dualis that is used when referring to two objects.

The empty set is obviously a set so it also needs a label and in English
it is {}``zero''. Like the other numerals the full understanding
of zero is far from integrated in most languages. For instance in
English singular is used when there is one object and plural when
there are more. Zero is not singular so it must be plural, i.e. one
man, two men, and zero men. There are other ways of indicating zero
and one can both say there was no man and there was no men which in
principle should be the same but in practice there is a light difference
in meaning. These are indications that the concepts of numerals and
of zero has not been fully incorporated in the language. For transfinite
cardinal numbers it is worse, but this is not the topic of this article.

The emergence of the concept of zero has a long development. It may
have started as an indication of an empty column on an abacus in ancient
China. When writing a number a blank seem to have been used to indicate
an empty column \cite{Needham1981,Martzloff1997}. Only later a symbol
for an empty space was used and the first evidence of using a symbol
for zero comes from India. Since then it has become an integrated
part of mathematics. The use of symbols better reflects modern standards
of mathematics than the languages that evolve quite slowly.

\section{Ordinal numbers}

When children start to learn numbers they do not distinguish between
cardinal numbers and ordinal number. The distinction only develops
gradually and some never really learn the difference. Still the difference
between cardinal numbers and ordinal numbers is far from well implemented
in the language. For instance in British English there is ambiguity
of writing 18 January or 18$^{th}$ January although in both cases
it should be clear that one refers to the 18$^{th}$ day of a month
called January rather than 18 month all called January. Sometimes
one also hear adults misusing ordinal number by saying that one competitor
perform double as bad as some other competitor in a competition because
one of the competitor became no. 6 and the other became no. 3. 

Let us take a look on the mathematical basis of ordinals. For two
well-ordered sets $A$ and $B$ there either exists an injective mapping
from $A$ to $B,$ or from $B$ to $A,$ such that the mapping maps
left sections into left sections. This mapping is unique so if we
have a large well-ordered set any element in any smaller well-ordered
set is mapped into a specific element in the large set. In order to
communicate an ordering we just have to communicate the mapping into
the larger set. It is convenient to have a fixed well-ordered set
and consider ordinal numbers as a well-ordered set of labels.

Therefore what people normally understand as ordinal numbers is a
well ordered class of labels that are used to label elements in a
well ordered set. We learn this as children. First we learn by heart
a sequence one, two, three,$\cdots$ and then we learn to count objects
in a set using this sequence. When counting we make a well-ordering
of the objects and assign labels from the sequence to the objects.
When we reach the last element in the set we automatically have the
label for the cardinality of the set. As counting is such an important
activity it is convenient that we reuse the cardinal numbers as ordinal
numbers but since cardinal numbers and ordinal numbers are used for
different purposes we modify the cardinal so that four become forth
and eleven become eleventh. Clearly the names of these labels are
derived from the names of the corresponding cardinal numbers except
for {}``first'' and {}``second''. This is an indicates that the
conceptual distinction between cardinal and ordinal number is a much
more recent conceptual invention than the cardinal numbers themselves.
Obviously there is no need for the label \textquotedbl{}zeroth\textquotedbl{}
in this system because an empty set has no element that should be
assigned a label like {}``zero'' or {}``zeroth''.

By experience we learn that counting leads to the same number independently
of the what well ordering was chosen. Only mathematicians would ever
realize that this could require a proof. For infinite sets the idea
of counting brakes down so what to do in this case. The first idea
might be simply to use the same labels for cardinal numbers and ordinal
number, but for infinite sets this is not possible. One problem is
that we do not know if the cardinal numbers are totally ordered. Only
if we assume the axiom of choice we can be sure that any two cardinals
are compatible. Another problem is that we do not know what the next
cardinal number is after $\aleph_{0}$ because this question is equivalent
to the continuum hypothesis. Therefore an attempt to use the cardinals
as labels would give serious restrictions to the foundation of mathematics
because adding the axiom of choice and the continuum hypothesis or
some versions of their negations have non-trivial consequences. These
problems show that cardinal and ordinal numbers are different types
of objects and that there is no reason to try to {}``reuse'' cardinal
numbers as ordinal numbers except if this is convenient.

It should be noted that people often use other labels to indicate
orderings. In sports, for instance, one use the well-ordered set \{gold,
silver, bronze\} for the three best competitors. Integers are used
to label other ordered sets as the set of Fourier coefficients $\dots c_{-2},c_{-1},c_{0},c_{1},c_{2},\dots$
for the Fourier series\[
f\left(x\right)=\sum_{n=-\infty}^{\infty}c_{n}\mathrm{e}^{inx}.\]
The theorem that for any two well-ordered sets there exists a unique
order isomorphism from a section of one of the sets to a left section
of the other set, cannot be extended to lattices or any other class
orderings. For instance there are many different order isomorphisms
between the ordered set \[
\left(\dots c_{-2},c_{-1},c_{0},c_{1},c_{2},\dots\right)\]
and $\mathbb{Z}.$ For special purposes there may be special reasons
to use a special enumeration. To start with $0$ is not more natural
than to start with $8$ or $-5.$ Sometimes there are even reasons
to use non-integers. For instance quantum particles associated with
irreducible representations of the rotation group are labeled by their
spin type that is one of the numbers $0,\nicefrac{1}{2},1,\nicefrac{3}{2},\cdots$
. Such special labels for a special purpose tells nothing about what
should be the default labeling. There is obviously nothing wrong about
using other sets than the ordinal numbers for labels but if we use
the ordinal numbers as labels we need to know exactly what the set
of ordinal numbers is. Words like {}``zeroth'' or {}``minus first''
that are sometimes used, but they indicate that the speaker has a
deviating opinion about what the ordinal numbers are. This just dilutes
the whole idea of having a fixed well-ordered set of elements that
can be used as labels for elements in other well-ordered sets.

\section{Order types}

There exists a preordering on well-ordered set where we define A to
be less than B if there exists a ordermorphish of A into B. The equivalence
classes are called order types. Cantor wanted labels for the order
types in parallel with the cardinal numbers that were as labels assigned
to equipotency classes. The ordering of order types is a well ordering
so one may reuse the labels for elements in a well-ordered set as
labels for order types. 

So here is the confusion. Ordinary people mainly use ordinal numbers
to label \emph{elements} of well-ordered sets but in enthusiasm of
his newly developed sets theory Cantor \cite{Cantor1897} introduced
ordinal numbers primarily as labels for \emph{order types.} For each
element in a well-ordered set $A$ one can assign a well ordered set
but there are two ways to do so. One is to use the set $\left\{ x\in A\mid x<a\right\} $
and another is to use $\left\{ x\in A\mid x\leq a\right\} .$ If the
second convention is used there is no element corresponding to the
empty set. If the first convention is used there is no element corresponding
to the whole set. Clearly the problem is that a finite well ordered
set with $n$ elements has $n+1$ left sections and there is no way
to completely solve this problem. Somehow the problem becomes kind
of invisible if we work with well ordered sets because a set can always
be embedded in a larger whereas one cannot find a subset of the empty
set. It was decided to associate the element $a$ with the set $\left\{ x\in A\mid x<a\right\} .$
Hence the first element was associated with the empty set and the
empty set has cardinality zero so among mathematicians zero has since
then been considered as the first ordinal number in complete contrast
with how people use ordinals to label elements of well ordered sets.
If we use the second convention instead the first section is the empty
set that has cardinality zero and the second section has cardinality
one etc. The first convention may appear more natural for infinite
sets but the correspondence between cardinals and ordinals breaks
done anyway when we reach infinite sets where sections with ordinal
types $\omega$ and $\omega+1$ both have cardinality $\aleph_{0}.$ 

The order types have a simple additive structure. One order type is
added to another order type by concatenating the corresponding well-ordered
sets. For infinite sets this gives a non-commutative addition. Although
addition makes sense for order types it makes no sense for elements
of well-ordered sets. If you are the third in one competition and
the fifth in a second competition you are not the seventh in any reasonable
sense although this is what Cantor ordinal numbers suggest. The way
to arrive at {}``seventh'' is by associating third with the well-ordered
set of the first and second in the first competition and then concatenate
this well-ordered set with well-ordered set consisting of the first,
the second, the third and the forth of the second competition. This
gives a well-ordered set of six elements and that is associated with
a seventh element. 

The real problem is that one wants to give labels to ordinal types
rather than to elements of well ordered sets as was the original idea
behind ordinal numbers. Actually I have never heard anybody outside
mathematical institutes use an ordinal number as a label for a set.
To do so would for instance be to use the word fifth about a ranking
list of four persons of a competition. Normally we use the words {}``before''
and {}``after'' rather than {}``less'' and {}``greater'' when
we refer to orderings like \textquotedbl{}the winner was Johnson and
after him came MacDonald\textquotedbl{}. I have never heard anybody
saying that a ranking list with three elements is \textquotedbl{}before\textquotedbl{}
a ranking list with five elements. 

To use the the ordinal numbers as labels for order types of sets is
clearly possible and perhaps relevant for special theoretical purposes
but it has very little to do with what people normally call ordinal
numbers.

\section{Construction and further confusion}

In the process of axiomatizing the foundation of mathematics the need
of constructing numbers emerged. The idea is that we should avoid
working with concepts that lead to inconsistencies. For instance after
defining a Hilbert space it is good practice to give at least one
example of a Hilbert space to ensure that all the consequences derived
from the axioms and definitions do not lead to contradictions. In
the attempts to axiomatize mathematics it was therefore needed to
specify a sequence of labels that could be used to label equipotent
sets. Using basic set theoretic constructions the following sequence
was derived:\begin{align*}
 & \varnothing\\
 & \left\{ \varnothing\right\} \\
 & \left\{ \varnothing,\left\{ \varnothing\right\} \right\} \\
 & \left\{ \varnothing,\left\{ \varnothing\right\} ,\left\{ \varnothing,\left\{ \varnothing\right\} \right\} \right\} \\
 & \vdots\end{align*}
 The idea was then to \emph{define} 3 as the label $\left\{ \varnothing,\left\{ \varnothing\right\} ,\left\{ \varnothing,\left\{ \varnothing\right\} \right\} \right\} .$
As proved by Gödel the use of such constructions is a not sufficient
to ensure consistency of mathematics. It is remarkable that the idea
has survived although it has completely failed and can only be used
to prove that there exists a set with say three elements if anybody
should doubt that. Although many textbooks still introduce the numbers
in this way most mathematicians would still count it as an error if
a studentwrote $2^{3}\subseteq3^{2}+2$ in an answer to an exercise.
The point is that cardinal numbers are primarily labels rather than
sets and a special language should be used for these labels (in this
case the use of $\leq$ rather than $\subseteq$).

John Conway introduced a different construction of numbers based on
combinatorial game theory \cite{Conway1976}. If playing games were
a more important activity than counting and ordering then his construction
should definitely be the basis of our notion of numbers. Conway constructs
something that he calls ordinal numbers in essentially the same way
as above. Nevertheless the additive structure on Conway ordinals is
different from the addition defined by concatenation of well-ordered
sets. For instance Conway's addition is commutative whereas addition
of ordinal number by concatenation is not. The most obvious conclusion
is that Conway ordinals are related to but not the same as what is
normally called ordinal numbers.

\section{Consequences}

The conclusion is that it is most natural to consider {}``zero''
as the first cardinal number and {}``first'' as the first ordinal
number. Mathematicians should simply not confuse order types with
ordinal numbers. Therefore mathematicians and computer scientists
should not by default enumerate like $x_{0},x_{1},x_{2},\ldots$ Similar
one should talk about a {}``polynomial of order zero'' rather than
a {}``zeroth order polynomial''.

This little story may also tell us why zero was invented much later
than one, two, and three. The understanding of zero requires an understanding
of the important distinction between cardinal numbers and ordinal
numbers. When counting we make an ordering of the set and assign labels
that are in principle ordinal numbers. When all elements have been
labeled we translate the last ordinal number into a cardinal number.
Hence zero was more difficult to understand than other cardinal numbers
because it did not correspond to any ordinal number.

\end{document}